\documentstyle[12pt]{article}

\newtheorem{fac0}{}[section]

\newtheorem{Def0}[fac0]{Definition}

\newtheorem{not0}[fac0]{Notation}

\newtheorem{propo2}{Proposition}[section]
\newtheorem{teo2}[propo2]{Theorem}

\newtheorem{ex2}[propo2]{Example}

\newtheorem{teo4}{Theorem}[section]
\newtheorem{rem4}[teo4]{Remark}

\newtheorem{lem4}[teo4]{Lemma}
\newtheorem{ex4}[teo4]{Example}
\newtheorem{tabla4}[teo4]{Table}

\newcommand{\qed}{\hfill\vbox{\hrule\hbox{\vrule\kern3pt
    \vbox{\kern6pt}\kern3pt\vrule}\hrule}} 

\font\Bbb=msbm10
\newcommand{\ZZ}{\mbox{\Bbb Z}}
\newcommand{\PP}{\mbox{\Bbb P}}

\newcommand{\OO}{\mbox{$O_{V_{d}}$}}

\begin{document}

\title{ VECTOR BUNDLES ON FANO 3-FOLDS WITHOUT INTERMEDIATE
COHOMOLOGY}
\author{
  $Enrique$  $ARRONDO$ \thanks{Partially supported by DGICYT PB96-0659.} \\
$AND$\\
$Laura$ $COSTA$ \thanks{Partially supported by DGICYT PB94-0850 and
 DGICYT PB96-0166.}
  \\ \\
 e-mail: enrique@sunal1.mat.ucm.es \\
 e-mail: costa@cerber.mat.ub.es    }
\date{DRAFT:5/4/98}
\maketitle

\section{Introduction}

A well-known result of Horrocks (see~ \cite{Ho}) states that a vector
bundle on a projective space has not intermediate cohomology if and
only if it decomposes as a direct sum of line bundles. There are two
possible generalizations of this result to arbitrary varieties. The
first one consists of giving a cohomological characterization of direct
sums of line bundles. This has been done for quadrics and Grassmannians
by Ottaviani (see ~\cite{Ot:1}, ~\cite{Ot:2}) and for rank-two vector
bundles on smooth hypersurfaces in $\PP^4$ by Madonna (see ~\cite{Ma}).

The second generalization, to which we make some contribution in this
paper, is to characterize vector bundles without intermediate
cohomology. Besides the result of Horrocks for projective spaces,
there is such a characterization for smooth quadrics due to Kn\"orrer
(see ~\cite{Kn}). In fact, it was shown by Buchweitz, Greuel and
Schreyer that only in the above two kind of varieties (projective
spaces and quadrics) there is, up to a twist by a line bundle, a finite
number of  indecomposable vector bundles without intermediate
cohomology (see ~\cite{BGS}). As a consequence, one should not expect to
find a precise characterization of vector bundles without intermediate
cohomology on arbitrary varieties. As far as we know, there is only a
result in that direction, for the Grassmann variety $G(1,4)$ of lines
in $\PP^4$, due to Gra\~na and the first author (see ~\cite{AG}).

\vspace{4mm}

The present paper deals with vector bundles without intermediate
cohomology on some Fano threefolds with the Betti number $b_2=1$.
Specifically, we will work over a smooth cubic in $\PP^4$, a smooth
complete intersection of type $(2,2)$ in $\PP^5$ and over a smooth
intersection of the Grassmannian $G(1,4)\subset \PP^9$ with three
hyperplanes. For rank-two vector bundles, we give a complete
classification (see Theorem~ \ref{rk2:sin}). For arbitrary rank, we
characterize which are the Chern classes of vector bundles without
intermediate cohomology and verifying some general conditions (see
Theorem~ \ref{superlopez}).

\vspace{4mm}

Next we outline the structure of the paper. In section 2, we
recall some general background. In section 3, we introduce some
standard rank-two vector bundles on our Fano threefolds without
intermediate cohomology, and we prove that any indecomposable rank-two
vector bundle is, up to a twist, one of these. In section 4, we give
the characterization of the Chern classes of ``sufficiently general''
vector bundles without intermediate cohomology and having arbitrary
rank.

\vspace{4mm}

\noindent{\bf Acknowledgements:} The authors want to express their
gratitude to the Dipartimento di Matematica of the Universit\`a degli
Studi di Milano, where most of the work has been developed. They also
want to stress the important help that the Maple package Schubert,
created by S. Katz and S.A. Str{\o}mme, has provided for some of the
computations needed in this work.

\section{Generalities.}
For $d=3,4,5$, let $V_d$ be a Fano 3-fold of degree $d$, index 2, with
the Betti number $b_2=1$ and such that one of the two generators of its
Picard group is spanned. Then, it has to be
\begin{itemize}
\item[] $V_3=$ a smooth cubic in $\PP^4$;
\item[] $V_4=$ a smooth complete intersection of type $(2,2)$ in
$\PP^5$;
\item[] $V_5=$ a smooth intersection of the Grassmannian $G(1,4)
\subset \PP^9$ with three hyperplanes.
\end{itemize}
Indeed, it is known (see e.g. \cite{class:1}) that there exist, up to a
deformation type, only four Fano threefolds of index two, the above ones
and a double covering $V_2$ of $\PP^3$ ramified along a quartic surface.

\vspace{4mm}

\begin{not0}
\rm
 Let $\OO(1)$ be the ample generator of $Pic(V_d)$. In our
situation, $\OO(1)$ is in fact very ample, and corresponds to the
hyperplane section $H$ of $V_d$ in $\PP^{d-1}$. Also, the canonical
divisor of $V_d$ is given by $K=-2H$. By a line on
$V_d$ we mean a rational curve $C$ such that $C\cdot \OO(1)=1$ and we
shall denote it by $L$, whereas a conic is a rational curve $C$ with
$C\cdot \OO(1)  =2$ and we shall denote it by $C$.

\vspace{4mm}
Since $H^2(V_d, \ZZ)$ is generated by the class of $H$,
$H^4(V_d, \ZZ)$ by a line $L$, and $H^6(V_d, \ZZ)$ by a point $P$,
we may identify the Chern classes $c_i \in H^{2i}(V_d, \ZZ)$
of vector bundles with integers, namely the coefficients of
$H$, $L$, and $P$. With this notation, $HL=1$, $H^2=dL$ and $H^3=d$.

Given $F$ a coherent sheaf on $V_d$ we will write $F(t)$ instead
of $F \otimes \OO(tH)$.
\end{not0}

\begin{Def0}
{\rm A vector bundle $F$ on $V_d$} has not intermediate cohomology {\rm
if and only if $H^{i}(V_d,F(t))=0$, for all $t \in \ZZ$ and $i=1,2$.}
 \end{Def0}

We will end this section by recalling the formulas for Chern classes of
a twist of a rank-$r$ vector bundle, the Riemann-Roch formula for
rank-$r$ vector bundles on $V_d$ and the well-known regularity
criterion of Castelnuovo-Mumford.

\begin{fac0} \label{Twist:Chern}
\rm
Given $F$ a rank-$r$ vector bundle on $V_d$ with Chern classes
$c_1(F)=c_1$, $c_2(F)=c_2$ and $c_3(F)=c_3$ it holds:
 \[ c_1(F(t))=c_1+rt; \]
  \[ c_2(F(t))=c_2+(r-1)tc_1d+{r\choose 2}t^2d; \]
  \[ c_3(F(t))=c_3+(r-2)c_2t+{r-1\choose 2}c_1t^2d+
     {r\choose 3}t^3d; \]
 \[ \chi(F)= \frac{dc_1^3-3c_1c_2}{6}+ \frac{dc_1^2-2c_2}{2}
+ \frac{(d+3)c_1}{3}+ \frac{c_3}{2}+r. \]
In particular,
\[ \chi(\OO(t))=\frac{dt^3}{6}+ \frac{dt^2}{2}+ \frac{(d+3)t}{3}+1.\]
\end{fac0}

\begin{fac0} \label{Castelnuovo:Mumford}
\rm
Let $O(1)$ be an ample invertible sheaf on a variety $X$ which is
generated by global sections. Let $F$ be a vector bundle with
\[ H^{i}(X,F(-i))=0 , \hspace{3mm} i=1,2, \cdots, dimX .\]
Then,
\begin{itemize}
\item[i)]
$ H^{i}(X,F(k))=0$, for any $ i=1,2, \cdots, dimX $ and $k \geq -
i$;
\item[ii)] $F$ is generated by global sections.
\end{itemize}
\end{fac0}

\section{Rank-two vector bundles without intermediate cohomology}

The aim of this section is to characterize, up to twist by line
bundles, rank-two vector bundles on $V_d$ without intermediate
cohomology.

In the next examples we will summarize  several facts about vector
bundles on $V_d$ that will be used in the sequel. Some of these results
are due to M. Szurek and J.A. Wi\'sniewski ([SW]).

\begin{ex2} \label{SL:resum} \rm Let $L$ be a line in $V_d$. Serre's
correspondence provides a vector bundle $S_L$ fitting in an exact
sequence:
 \refstepcounter{equation}
\[ 0 \rightarrow \OO
\rightarrow   S_L  \rightarrow
 I_L  \rightarrow
 0  \hspace{8mm} (\theequation) \label{se3}. \]
The rank-two vector bundle $S_L$ has Chern classes $(c_1,c_2)=(0,1)$, is
semistable and has not intermediate cohomology. It holds that $S_L$ has
only one section, whereas $S_L(1)$ is generated by global sections.
The latter comes from the exact sequence (\ref{se3}), the fact that
$\OO(1)$ and $I_L(1)$ are generated by its sections and the vanishing of
$H^1(\OO(1))$.
\end{ex2}

\begin{ex2} \label{SC:resum} \rm If $C$ is a conic in $V_d$, let $S_C$
be the rank-two vector bundle obtained from Serre's construction and
fitting in a non-trivial extension:
 \refstepcounter{equation}
\[ 0 \rightarrow \OO(-1)
\rightarrow   S_C  \rightarrow
 I_C  \rightarrow
 0  \hspace{8mm} (\theequation) \label{se1}. \]
It holds that $S_C$ is stable and generated by its global sections. It
has not intermediate cohomology and has Chern classes
$(c_1,c_2)=(-1,2)$.

When $d=4$, $S_C$ is related to the spinor bundles of the quadrics
containing $V_4$. More precisely, let $\Pi$ be the unique plane
containing $C$. Then, only one quadric $Q_4$ containing $V_4$ contains
$\Pi$. If $Q_4$ is smooth, then $S_C(1)$ is the restriction to $V_d$ of
one of the spinor bundles of $Q_4$. In other words, if we identify
$Q_4$ with the Grassmann variety $G(1,3)$ of lines in $\PP^3$, then
$S_C(1)$ is the restriction to $V_4$ of one of the universal bundles of
$G(1,3)$. On the other hand, when $d=5$, $S_C(1)$ is the restriction to
$V_5$ of the rank-two universal quotient bundle on $G(1,4)$.
\end{ex2}

\begin{ex2} \label{SE:resum} \rm Using again Serre's correspondence, an
elliptic curve $E$ of degree $d+2$ gives rise to a stable, rank-two
vector bundle $S_E$, with Chern classes $(c_1,c_2)=(0,2)$, without
intermediate cohomology and appearing in an exact sequence:
 \refstepcounter{equation}
\[ 0 \rightarrow \OO(-1)
\rightarrow   S_E  \rightarrow
 I_E(1)  \rightarrow
 0  \hspace{8mm} (\theequation) \label{se2}. \]
It also holds that $S_E(1)$ is generated by its global sections.
\end{ex2}

\vspace{4mm}
Next we will prove that the above examples are essentially the only
rank-two vector bundles on $V_d$ without intermediate cohomology.

\begin{teo2}
\label{rk2:sin}
An indecomposable rank-two vector bundle $F$ on $V_d$ with $d=3,4,5$ has
not intermediate cohomology if and only if it is a twist of either
$S_L$, or $S_C$, or $S_E$.
\end{teo2}
\noindent
{\em Proof:}
Let $F$ be an indecomposable rank-two vector bundle on $V_d$
with Chern classes $c_1(F)=c_1$ and $c_2(F)=c_2$, which has not
intermediate cohomology. Without loss of generality we can assume that
$F(-1)$ has no sections while $F$ has. Hence, any section of $F$
gives rise to an exact sequence:
 \refstepcounter{equation}
\[ 0 \rightarrow \OO
\rightarrow  F \rightarrow
I_{D}(c_1) \rightarrow
 0  \hspace{8mm} (\theequation) \label{se10} \]
where $D$ is a scheme of pure dimension one. Observe that it cannot be
$D=\emptyset$, since Ext$^1(\OO(c_1),\OO)\cong H^1(\OO(-c_1))=0$, and
this would imply that $F$ decomposes.
\bigskip

\noindent
{\em Claim:}
\begin{center}
$F(2-c_1)$ is generated by its global sections and therefore $c_1 \leq
2.$
\end{center}
\noindent
{\em Proof of the Claim:} Using Serre's duality and the natural
identification $F^*\cong F(-c_1)$ we have:
\[ h^3(F(-1-c_1))=h^0(F(-1))=0 .\]
Hence, since $F$ has not intermediate cohomology, using
Castelnuovo-Mumford criterion~ \ref{Castelnuovo:Mumford}, we obtain the
claim.

\vspace{4mm}
 We will prove the theorem distinguing
two cases. \newline

\vspace{4mm}
\noindent
{\em Case A:} Assume that $c_1 \leq 0$.
{}From Serre's duality and (\ref{se10}) twisted with $\OO(-c_1-1)$, we
get that the only cohomology of $F(-1)$ is
$h^3(F(-1))=h^0(F(-c_1-1))=h^0(\OO(-c_1-1))$. Computing $\chi(F(-1))$
and $\chi(\OO(-c_1-1))$
from \ref{Twist:Chern} we get:
\[ \chi(F(-1))=\frac{c_1^3d}{6}-\frac{c_1c_2}{2}-\frac{dc_1}{6}+c_1 ;\]
\[ \chi(\OO(-c_1-1))=-\frac{c_1^3d}{6}+\frac{dc_1}{6}-c_1 .\]
Hence, from $\chi(\OO(-c_1-1))=-\chi(F(-1))$
we obtain the identity $c_1c_2=0$. Since
$D\neq\emptyset$, we have $c_2\neq 0$, and hence $c_1=0$. Again from
\ref{Twist:Chern} and the fact that $F$ has not intermediate cohomology
we get $2-c_2=\chi(F)=h^0(F)\ge 1$, from which $c_2=1$. This means that
$D$ is a line $L$, and therefore $F \cong S_L$.

\vspace{4mm}
\noindent
{\em Case B:} If $c_1>0$, we know from the claim that the only
possibilities are $c_1=1,2$. In fact, this can be obtained, together
with a much more precise information, with the same techniques as in
case A. Indeed, we have now that $h^0(F(-c_1))=0$, so that
$h^3(F(-1))=h^3(F(-2))=0$. Therefore, $\chi(F(-1))=\chi(F(-2))=0$.
Again, using \ref{Twist:Chern}, we have
\[ \chi(F(-2))=\frac{c_1^3d}{6}-\frac{c_1c_2}{2}+c_1
               -\frac{dc_1^2}{2}+c_2-2 +\frac{dc_1}{3}. \]
Comparing with the above expression for $\chi(F(-1))$, we obtain that
$c_2-\frac{dc_1^2}{2}-2 +\frac{dc_1}{2}=0$, and therefore the only
possibilities are
$(c_1,c_2)=(1,2),(2,d+2)$.

\vspace{4mm}
In the first case, $D$ has degree $c_2=2$, and by Serre's correspondence
$\omega_D\cong O_D(c_1-2)=O_D(-1)$, which means that $D$ is a
plane conic $C$ (maybe singular or even non-reduced). Therefore, $F
\cong S_C(1)$.

In the second case, by the claim, $F$ is then generated by its global
sections. This implies that a general section of $F$ vanishes on a
smooth curve $D$. Since $\omega_D=O_D(c_1-2)=O_D$, this means that $D$
is an elliptic curve. Its degree is $c_2=d+2$. Therefore,
$F \cong S_E(1)$ where $S_E$ is the vector bundle on $V_d$ associated
to an elliptic curve $E$ of degree $d+2$.  \qed

\section{Vector bundles of higher rank without intermediate cohomology}

In this section we will study vector bundles of rank $r\geq 3$ without
intermediate cohomology. While in the previous section we gave a
complete classification, in this section we will just characterize the
possible Chern classes. We will also need to make some general
assumptions.

\vspace{4mm}

\noindent{\bf Convention:} For simplicity, we will
frequently use $(\star)$ to denote the following conditions on a
rank-$r$ vector bundle $F$ on $V_d$:
\begin{itemize}
\item[$(\star)$] The vector bundle $F$ has not trivial summands,
$h^0(F)\geq r$, $h^0(F(-1))=0$ and $F$ has $r-1$ sections whose
dependency locus has codimension two.
\end{itemize}

\begin{lem4}\label{suma-star} Let $F'$ and $F''$ be two vector bundles
on $V_d$ with respective ranks $r'$ and $r''$. If $F'$ and $F''$ verify
$(\star)$, then also $F'\oplus F''$ verifies $(\star)$.
\end{lem4}

\noindent {\em Proof:} It is clear from the definition that, since
$F'$ and $F''$ verify $(\star)$, then their global sections generate
them up to probably a subvariety of codimension at least two.
Therefore, the same holds for $F'\oplus F''$. This implies that the
dependency locus of $r'+r''-1$ general sections of $F'\oplus F''$ has
codimension two. {}From this, it follows immediately that $F'\oplus F''$
verifies $(\star)$. \qed

\vspace{4mm}

Before stating and proving our main theorem, we will give a series of
examples. Eventually this list will become the complete list for small
rank.

\begin{ex4}\label{c1=1} \rm Let $D$ be a rational curve of degree $r\leq
d+1$ in $V_d$. Serre's construction yields a rank-$r$ vector bundle $F$
fitting in an exact sequence
\[0\to O_{V_d}^{r-1}\to F\to I_D(1)\to 0.\]
Since $D$ is projectively normal in $\PP^{d+1}$ if and only if $r\leq
d+1$, then $F$ has not intermediate cohomology. From the above exact
sequence it holds that $c_1(F)=1$, $c_2(F)=r$ and $c_3(F)=r-2$, and
also $h^0(F)=d$. Notice that, when $d=5$ and $r=3$, $F$ is the
restriction to $V_5$ of the dual of the rank-three universal subbundle
on $G(1,4)$, and in particular it is generated by its global sections.
\end{ex4}

\begin{ex4}\label{c1=2,r=3} \rm Take $F$ to be the rank-three vector
bundle given by the above example when $r=3$. Then $\bigwedge^2F\cong
F^*(1)$ has not intermediate cohomology, has general sections and
Chern classes $c_1=2$, $c_2=d+3$ and $c_3=2$. Since, for $d=5$, $F$
is generated by its global sections, then the same holds for
$\bigwedge^2F$.
\end{ex4}

\begin{ex4}\label{c1=3,r=3} \rm Consider the elliptic curve $E$
given in Example~ \ref{SE:resum}. Exact sequence (\ref{se2}) implies that
there is a section $s\in H^0(S_L(2))$ vanishing on $E$ and not vanishing
on a divisor of $V_d$. If $D'$ is the residual curve of $E$ inside the
zero locus of $s$, it turns out that $D'$ has degree $3d-1$ and
(arithmetic) genus $2d-2$. The standard mapping cone construction (see
for instance Lemma 3.2 in~ \cite{tesis})
provides an exact sequence
\[0\to S_E(-3)\to\OO(-2)\oplus S_L(-2)\to I_{D'}\to 0.\]
We use now this exact sequence to see that there is a section $s'\in
H^0(S_C(3))$ vanishing on $D'$ and not vanishing on a divisor of $V_d$.
We have now the followig exact sequences for $D$, the residual curve of
$D'$ inside the zero locus of $s'$:
 \refstepcounter{equation}
\[0\to S_L(-3)\oplus\OO(-3)\to S_E(-2)\oplus S_C(-2)\to I_D\to 0
\hspace{8mm} (\theequation) \label{se5}\]
and its dual
 \refstepcounter{equation}
\[0\to \OO\to S_E(-2)\oplus S_C(3)\to S_L(3)\oplus\OO(3)\to
\omega_D\to 0. \hspace{8mm} (\theequation) \label{se6}\]
The curve $D$ has degree $3d+3$ and arithmetic genus $2d+4$. {}From the
exact sequence (\ref{se6}) tensored with $\OO(-3)$ it follows that
$h^0(\omega_D(-1))=2$, so that Serre's construction yields a rank-three
vector bundle $F$ with Chern classes $c_1=3$, $c_2=3d+3$, $c_3=d+3$ and
fitting in an exact sequence
 \refstepcounter{equation}
\[0\to O_{V_d}^2\to F\to I_D(3)\to 0. \hspace{8mm} (\theequation)
\label{se7} \] Comparing with (\ref{se5}) it easily follows that
$H^0(F(-1))=0$,
$H^1(F(l))=0$ for any $l\in{\Bbb Z}$, and $H^2(F(l))=0$ for any $l\geq
-1$. It also follows from Serre's construction that $H^2(F(-2))=0$.
Hence, Mumford-Castelnuovo criterion applied to $F^*(2)$ will imply
that $F$ has not intermediate cohomology as soon as we prove that
$h^1(F^*(1))=h^2(F(-3))=0$. We observe that, from the exact sequences
(\ref{se6}) and the dual of (\ref{se7}) and our construction, there is
a commutative diagram
\[\matrix{
0\to&\OO(-2)&\to&F^*(1)    &\to&\OO(1)^2  &\to&\omega_D&\to 0\cr
    &  ||   &   &\downarrow&   &\downarrow&   &   ||   &     \cr
0\to&\OO(-2)&\to&S_E\oplus
S_C(1)&\to&S_L(1)\oplus\OO(1)&\to&\omega_D&\to 0}\]
The vanishing we want to prove reduces to the surjectivity of the map
$H^0(\OO(1)^2)\to H^0(\omega_D)$, but this is proved by just taking
cohomology in the above diagram.

Notice that the above construction does not imply a priori that the
curve $D$ is smooth. However, once the vector bundle $F$ is
constructed, it follows by the Castelnuovo-Mumford criterion that $F$
is generated by its global sections. Hence, the dependency locus of
two general sections of $F$ is a smooth curve $D$.

When $d=5$, such a vector bundle $F$ can be constructed in a more
direct way. Indeed, in this case $S_C(1)$ is the restriction to $V_5$ of
the universal rank-two quotient bundle $Q$ on $G(1,4)$. Since $S^2Q$ has
not intermediate cohomology, neither $F=S^2(S_C(1))$ has. It follows at
once that $F$ is a rank-three vector bundle, generated by its global
sections and its Chern classes are $c_1=3$, $c_2=18$ and $c_3=8$.
\end{ex4}


\begin{ex4}\label{c1=2,r=7,d=5} \rm Let $G$ be the rank-three vector
bundle obtained in Example~\ref{c1=2,r=3} when $d=5$. Since $G$ was
generated by its global sections and $h^0(G)=10$, then there is an exact
sequence
\[0\to K\to O_{V_5}^{10}\to G\to 0\]
where $K$ is defined as a kernel. It follows that $F=K^*$ is a
rank-seven vector bundle without intermediate cohomology and generated
by its global sections. It also holds that $c_1(F)=2$, $c_2(F)=12$ and
$c_3(F)=10$.
\end{ex4}

\newpage

\begin{tabla4}\label{tabla} \rm From Lemma~ \ref{suma-star} and the
above examples, we can give the following list of rank-$r$ vector
bundles $F_{r,c}$ without intermediate cohomology, verifying
$(\star)$ and with $c_1(F_{r,c})=c$

\begin{table}[hbt]
\begin{center}
\begin{tabular}{|c|c|c|c|}
\hline {$d$}&{$r$}&$(c,c_2,c_3)$&$F_{r,c}$\\
\hline
\hline
3,4,5&3&$(1,3,1)$& Example~\ref{c1=1}\\
\hline
3,4,5&3&$(2,d+3,2)$& Example~\ref{c1=2,r=3}\\
\hline
3,4,5&3&$(3,3d+3,d+3)$& Example~\ref{c1=3,r=3}\\
\hline
  4,5&4&$(1,4,2)$& Example~\ref{c1=1}\\
\hline
3,4,5&4&$(2,d+4,4)$& $S_C(1)\oplus S_C(1)$\\
\hline
3,4,5&4&$(3,3d+4,d+6)$& $S_C(1)\oplus S_E(1)$\\
\hline
3,4,5&4&$(4,6d+4,4d+8)$& $S_E(1)\oplus S_E(1)$\\
\hline
    5&5&$(1,5,3)$& Example~\ref{c1=1}\\
\hline
3,4,5&5&$(2,d+5,6)$& $S_C(1)\oplus F_{3,1}$\\
\hline
3,4,5&5&$(3,3d+5,d+9)$& $S_C(1)\oplus F_{3,2}$\\
\hline
3,4,5&5&$(4,6d+5,4d+2)$& $S_C(1)\oplus F_{3,3}$\\
\hline
3,4,5&5&$(5,10d+5,10d+15)$& $S_E(1)\oplus F_{3,3}$\\
\hline
3,4,5&6&$(2,d+6,8)$& $F_{3,1}\oplus F_{3,1}$\\
\hline
3,4,5&6&$(3,3d+6,d+12)$& $S_C(1)\oplus S_C(1)\oplus S_C(1)$\\
\hline
3,4,5&6&$(4,6d+6,4d+16)$& $S_C(1)\oplus S_C(1)\oplus S_E(1)$\\
\hline
3,4,5&6&$(5,10d+6,10d+20)$& $S_C(1)\oplus S_E(1)\oplus S_E(1)$\\
\hline
3,4,5&6&$(6,15d+6,20d+24)$& $S_E(1)\oplus S_E(1)\oplus S_E(1)$\\
\hline
  4,5&7&$(2,d+7,10)$& $F_{4,1}\oplus F_{3,1}$\\
\hline
3,4,5&7&$(3,3d+7,d+15)$& $S_C(1)\oplus S_C(1)\oplus F_{3,1}$\\
\hline
3,4,5&7&$(4,6d+7,4d+20)$& $S_C(1)\oplus S_E(1)\oplus F_{3,1}$\\
\hline
3,4,5&7&$(5,10d+7,10d+25)$& $S_E(1)\oplus S_E(1)\oplus F_{3,1}$\\
\hline
3,4,5&7&$(6,15d+7,20d+30)$& $S_E(1)\oplus S_E(1)\oplus F_{3,2}$\\
\hline
3,4,5&7&$(7,21d+7,35d+35)$& $S_E(1)\oplus S_E(1)\oplus F_{3,3}$\\
\hline
\end{tabular}
\end{center}
\end{table}

\end{tabla4}

\vspace{4mm}

\begin{rem4} \rm Notice that in the above table, vector bundles with the
same invariants can be constructed in different ways. Moreover, in
some cases it is possible to consider vector bundles given by
non-trivial extensions instead of direct sums. For instance, it is
immediate to see that $Ext^1(S_E(1),S_E(1))\neq 0$. In fact, the
dimension of this space of extensions is five, the same as the
dimension of the space of elliptic curves of degree $d+2$ contained in
$V_d$. Clearly, a vector bundle obtained from a non-trivial extension
also verifies condition $(\star)$ and has not intermediate cohomology.
\end{rem4}

\vspace{4mm}

\begin{rem4} \rm
In case of rank $r=3$, it is clear that the vector bundles in the
above list are indecomposable, since otherwise any direct summand of
rank one would be necessarely trivial. In case $r=4$ and $c_1=1$
(even when $d=3$, as constructed in Example~ \ref{c1=1}) the
corresponding vector bundle in the above list is also indecomposable.
Indeed, we can rule out (using the classification given in the
previous section) all its possible rank-two direct summands. What seems
to be difficult is to decide in general when it is possible to find an
indecomposable vector bundle without intermediate cohomology and having
some of the above invariants. This is so even when we know that we can
take non-trivial extensions as in the above remark.
\end{rem4}

\vspace{4mm}
We can now state and prove the main result of this section.

\begin{teo4} \label{superlopez} Let $r\geq 3$ be an integer. There
exists a rank-$r$ vector bundle $F$ on $V_d$ without intermediate
cohomology, verifying $(\star)$ and with Chern classes $c_1(F)=c_1$,
$c_2(F)=c_2$ and
$c_3(F)=c_3$ if and only if
\[ (i)\ c_2=\frac{dc_1^2}{2}+r-\frac{dc_1}{2};\]
\[ (ii)\
c_3=-2c_1+c_1r-\frac{dc_1^2}{2}+\frac{dc_1^3}{6}+\frac{dc_1}{3};\]
\[ (iii)\ \frac{r}{d}\leq c_1\leq r.\]
Moreover, in this situation, the dependency locus of $r-1$ general
sections of such $F$ is a curve $D$ of degree $c_2$ and arithmetic genus
\[p_a(D)=(c_1-1)(r-1)-dc_1^2+\frac{dc_1^3}{3}+\frac{2dc_1}{3}. \]
\end{teo4}

\noindent
{\em Proof:}
Let $F$ be a rank-$r$ vector bundle on $V_d$ verifying $(\star)$, with
Chern classes $c_1(F)=c_1$, $c_2(F)=c_2$ and $c_3(F)=c_3$, which has
not intermediate cohomology. By assumption, $F$ has $r-1$ sections
whose dependency locus has codimension two. These sections of $F$ give
rise to an exact sequence:
 \refstepcounter{equation}
\[ 0 \rightarrow \OO^{r-1}
\rightarrow  F \rightarrow
I_{D}(c_1) \rightarrow
 0  \hspace{8mm} (\theequation) \label{se11} \]
where $D$ is a scheme of pure dimension one. Observe that it cannot be
$D=\emptyset$, since Ext$^1(\OO(c_1),\OO^{r-1})\cong H^1(\OO(-c_1))
^{\oplus(r-1)}=0$, and this would imply that (\ref{se11}) splits, and
hence $F$ would have a trivial summand.

\vspace{4mm}
Let us see that $c_1>0$. Indeed, assume that $c_1 \leq 0$.
{}From Serre's duality and (\ref{se11}) twisted by $\OO(-1)$,
the only cohomology of $F(-1)$ is
$h^3(F(-1))=h^3(\OO(c_1-1))$. Since $c_1(F(-1))=c_1-r$,
$c_2(F(-1))=c_2-(r-1)c_1d+{r \choose 2}d$ and
$c_3(F(-1))=c_3-(r-2)c_2+{r-1 \choose 2}c_1d- {r \choose 3}d$,
 computing $\chi(F(-1))$ and $\chi(\OO(c_1-1))$
from \ref{Twist:Chern} we get:
\[ \chi(F(-1))=\frac{c_1^3d}{6}-\frac{c_1c_2}{2}-\frac{dc_1}{6}+c_1
+ \frac{c_3}{2} ;\]
\[ \chi(\OO(c_1-1))=\frac{c_1^3d}{6}-\frac{dc_1}{6}+c_1 .\]
Hence, from $\chi(\OO(c_1-1))=\chi(F(-1))$
we obtain the identity $ c_3-c_1c_2=0$.

Similarly, we have $\chi(F(-c_1)) =(r-1)\chi(\OO(-c_1))$, which is
equivalent to:
$$
(-\frac{(r-1)c_1^3}{6}-\frac{(r-1)c_1}{3}+\frac{(r-1)c_1^2}{2})d
+\frac{c_1c_2}{2}-c_2+\frac{c_3}{2}-rc_1+c_1+r=$$
$$ (r-1)((-\frac{c_1^3}{6}+\frac{c_1^2}{2}
-\frac{c_1}{3})d-c_1+1)$$
and $\chi(F(-c_1-1))=(r-1)\chi(\OO(-c_1-1)$, which is equivalent to:
$$
(-\frac{(r-1)c_1^3}{6}+\frac{(r-1)c_1}{6})d
+\frac{c_1c_2}{2}+\frac{c_3}{2}-rc_1+c_1=$$
$$(r-1)((-\frac{c_1^3}{6}
+\frac{c_1}{6} )d-c_1)$$
from which we respectively obtain the identities $c_1c_2+c_3=c_2-1$
and $c_1c_2+c_3=0$. Putting together the three identities, we get
$c_1=c_3=0$ and $c_2=1$, so that in particular $D$ is a line $L$.

\vspace{4mm}
The dual of (\ref{se11}) is
 \refstepcounter{equation}
\[ 0 \rightarrow \OO(-c_1)
\rightarrow  F^* \rightarrow\OO^{\oplus(r-1)}\to
\omega_{D}(2-c_1) \rightarrow
 0.  \hspace{8mm} (\theequation) \label{se12} \]

Since $\omega_{D}(2-c_1)$ is trivial in our case, it follows that
$F^*$ has $r-2$ independent sections such that the composed map
$\OO^{\oplus(r-2)}\to F^*\to \OO^{\oplus(r-1)}$ defines a direct
summand. But this implies that $F\cong \OO^{\oplus(r-2)}\oplus S_L$,
which is a contradiction as $r\ge 3$. Hence, $c_1>0$.

\vspace{4mm}
Since $c_1>0$, we have that $h^3(F(-1))=h^3(F(-2))=0$ (the latter
equality is equivalent to $h^0(F^*)=0$, which is obtained from
(\ref{se12}) as above, by using that $F$ has no trivial summands).
Therefore, $\chi(F(-1))=\chi(F(-2))=0$. Again, using
\ref{Twist:Chern}, we have:
\[ \chi(F(-1))=\frac{c_1^3d}{6}-\frac{c_1c_2}{2}-\frac{dc_1}{6}
+c_1 +\frac{c_3}{2};\]
\[ \chi(F(-2))=\frac{c_1^3d}{6}-\frac{c_1c_2}{2}+c_1
               -\frac{dc_1^2}{2}+c_2-r +\frac{dc_1}{3}+ \frac{c_3}{2}.
\]
The vanishing of the above two terms immediately implies that
$c_2(F)$ and $c_3(F)$ verify the identities $(i)$ and $(ii)$.

Using once more \ref{Twist:Chern}, we obtain that
$-h^3(F(-3))=\chi(F(-3))=d(c_1-r)$. This implies that $c_1\le r$.
The other inequality in $(iii)$ follows from the fact that $h^0(F)\geq
r-1$, Riemann-Roch equality~ \ref{Twist:Chern} and the identities $(i)$
and $(ii)$.

Finally, the genus of the curve $D$ is $\chi(I_D)$, which can be
computed from (\ref{se11}), and it is
\[(c_1-1)(r-1)-dc_1^2+\frac{dc_1^3}{3}+\frac{2dc_1}{3}. \]

\vspace{4mm}

Conversely, let $r$, $c_1$, $c_2$ and $c_3$ be integers verifying
$(i)$, $(ii)$ and $(iii)$. We will prove by induction on $r$ that
there exists a rank-$r$ vector bundle $F$ on $V_d$ without
intermediate cohomology, verifying $(\star)$ and with Chern classes
$c_1(F)=c_1$, $c_2(F)=c_2$ and $c_3(F)=c_3$. For $r\leq 7$ this is
showed in Table~ \ref{tabla}, so let us assume $r\geq 8$. In this
case, the proof goes as for $r\leq7$, using Lemma~ \ref{suma-star}. We
will distinguish three cases:
\bigskip

\noindent {\em Case A)} If $c_1=r$, we can take $F=S_E(1)\oplus F'$,
where $F'$ is a rank-$(r-2)$ vector bundle on $V_d$ without
intermediate cohomology, verifying $(\star)$ and with first Chern class
$c_1(F')=r-2$. Indeed, such $F'$ can be constructed by induction
hypothesis.
\bigskip

\noindent {\em Case B)} If $\frac{r-2}{d}+1\leq c_1\leq r-1$, we
take $F=S_C(1)\oplus F'$, where $F'$ is a rank-$(r-2)$ vector bundle on
$V_d$ without intermediate cohomology, verifying $(\star)$ and with first
Chern class $c_1(F')=c_1-1$. This vector bundle $F'$ exists by
induction hypothesis, since $F'$ verifies $(iii)$.
\bigskip

\noindent {\em Case C)} If $\frac{r}{d}\leq c_1\leq
\frac{r-1}{d}+1$, we take now $F=F'\oplus F''$, where $F'$ and $F''$
are vector bundles on $V_d$ of respective ranks $r-d$ and $d$ without
intermediate cohomology, verifying $(\star)$ and with Chern classes
$c_1(F')=c_1-1$, $c_1(F'')=1$. Once again, $F'$ verifies the
inductive assumption since $r\geq 8$, and $F''$ was constructed
in Example~ \ref{c1=1}.

\qed

\vspace{4mm}

\begin{rem4} \rm If we do not make the assumption that our vector
bundles have at least $r$ independent sections, the same proof is
clearly valid to see that in that case $(i)$ and $(ii)$ holds, and
instead of $(iii)$ we would have that $\frac{r-1}{d}\leq c_1\leq r$.
What presents some difficulty is to prove the actual existence of
vector bundles with those invariants. This is due to the fact that we
cannot make use of Lemma~ \ref{suma-star}. In fact, we are still able
to prove the existence of most of them, but some cases still remain
unknown to us.
\end{rem4}

\vspace{4mm}

\begin{rem4} \rm Our results are in fact valid over any Fano threefold
of degree $d$ and index two and arbitrary Betti number. In that case,
the only assumption we need on the vector bundles is their determinant
to be a multiple of $\OO(1)$, the hyperplane section of $V_d$ in
$\PP^{d-1}$.
\end{rem4}

\end{document}